\documentclass[12pt]{conm-p-l}
\usepackage{amsthm,amssymb,amsmath,calc,eucal,ifthen}
\usepackage[matrix,arrow]{xy}

\newtheorem{theorem}{Theorem}[section]
\newtheorem{lemma}[theorem]{Lemma}
\newtheorem{proposition}[theorem]{Proposition}
\newtheorem{corollary}[theorem]{Corollary}

\theoremstyle{remark}
\newtheorem{remark}[theorem]{Remark}

\newtheoremstyle{rmdefinition}{}{}{\upshape}{}{\bfseries}{.}{ }{}

\theoremstyle{rmdefinition}

\newcommand{\derived}[2][1]{\ifthenelse{\equal{#1}{1}}{{#2}^{\prime}}{\ifthenelse{\equal{#1}{2}}{{#2}^{\prime\prime}}{\ifthenelse{\equal{#1}{3}}{{#2}^{\prime\prime\prime}}{{#2}^{(#1)}}}}}

\DeclareMathOperator{\tr}{tr}

\DeclareMathOperator{\spec}{spec}

\newcommand{\dml}{\Delta_\sigma}

\begin{document}
\title[Asymptotic spectral 
properties for discrete magnetic Laplacians]{Kato's inequality and asymptotic spectral 
properties for discrete magnetic Laplacians}

\author{J\'ozef Dodziuk}
\address{Ph.D.\ Program in Mathematics, Graduate Center of CUNY, New York, NY 10016, USA.}
\email{jozek@derham.math.qc.edu}
\author{Varghese Mathai}
\address{
Department of Mathematics, University of Adelaide, Adelaide 5005,
Australia}
\email{vmathai@maths.adelaide.edu.au}

\begin{abstract}In this paper, a discrete form of the Kato inequality
for discrete magnetic Laplacians on graphs is used to study asymptotic
properties of the spectrum of discrete magnetic Schr\"odinger operators.
We use the existence of a ground
state with suitable properties for the ordinary combinatorial Laplacian
and semigroup domination to relate the combinatorial Laplacian with
the discrete magnetic Laplacian. Our techniques yield existence and
uniqueness of the fundamental solution for the heat kernel on a graph of
bounded valence.

\end{abstract}

\maketitle

\section*{Introduction}

The discrete magnetic Laplacian (DML), $\Delta_\sigma$ 
and discrete magnetic Schr\"odinger operator (DMSO), $\Delta_\sigma + T$, 
where $T$ is in the commutant of a projective $G$ action,
on graphs had been studied for many years by condensed matter and solid state
physicists, mainly on graphs having a free action of a crystallographic group $G$
cf. \cite{Bellissard}. In particular,
these operators are Hamiltonians for the discrete model of the quantum Hall effect. Sunada
\cite{Sunada}
extended the definition of these operators to arbitrary graphs equipped 
with a free action of a general finitely generated discrete group $G$.

In this paper, we establish asymptotic properties of the spectrum of DMLs and 
DMSOs,
such as the large time decay of the $G$-trace of the heat kernels of these operators,
and also  the positivity of the Fuglede-Kadison determinant of these operators - 
these results could have  applications to physics, particularly to 
analysis of the fractional quantum Hall effect. 

A discrete form of the Kato inequality is used in the paper. Other
applications of it are in \cite{Colin}, for finite graphs 
and by physicists in the case of graphs having a free action of a crystallographic group.
Together with semigroup domination techniques and some well known results of Varopoulos
\cite{VCC}, the desired asymptotic properties of the spectrum of DMLs and DMSOs of the form 
$\Delta_\sigma + T$, where $T\ge 0$  is established here. This result 
is optimal in the case when $G$ is amenable - here the characterization of amenability 
of $G$ via  the bottom of the spectrum of the discrete Laplacian $\lambda_0$ 
due to R. Brooks \cite{Brooks1982}, \cite{Sunda-Polly}
is used. In the nonamenable case, we establish basic properties of the 
ground state of the discrete Laplacian, and use it to establish 
asymptotic properties of the spectrum of DMSOs of the form 
$\Delta_\sigma + T$, where $T\ge -\lambda_0$.

\section{Kato's inequality}

We show that the technique of domination of semigroups as developed in
\cite{HSU1977} can  be used to relate discrete magnetic Laplacians on a
graph with the ordinary discrete Laplacian. Our proofs follow
closely \cite{Berard1986}. We begin with the positivity preserving
property of the ordinary discrete Laplacian $$
\Delta : C^0_{(2)}(K) \longrightarrow C^0_{(2)}(K)$$
acting on the space of square-integrable cochains on the graph $K$,
i.e.\ the space of complex-valued functions $f$ on vertices of $K$ such 
that
$\sum_v |f(v)|^2 < \infty$.
For a function $f\in C^0_{(2)}(K)$ and a vertex $v$ of $K$
\begin{equation}\label{lapl}
\Delta f(v) = \sum_{w \sim v} (f(v) - f(w))
\end{equation}
where the notation $w\sim v$ denotes that $v$ and $w$ are connected by
an edge.

\begin{lemma} \label{positivity}
Suppose $f \in C^0_{(2)}(K)$ is a real-valued function
and $\lambda >0$. If $(\Delta + \lambda I)f
\geq 0$, then $f\geq 0$ i.e.\ $(\Delta +\lambda I)^{-1}$ is positivity
preserving.
\end{lemma}
\begin{proof}
Suppose $f$ satisfies the assumptions of the lemma but attains a
negative value. Since $f \in C^0_{(2)}(K)$ it vanishes at infinity and
therefore attains a negative minimum at a vertex $v_0$ of $K$. Now, by
(\ref{lapl}) $$
(\Delta + \lambda I)f(v_0) = \sum_{w \sim v_0} \left (f(v_0) - 
f(w)\right
)
+\lambda f(v_0).$$
The left-hand side of this equality is nonnegative by assumption whereas
the right-hand side is strictly negative. The contradiction proves the
lemma.
\end{proof}

Now let $\Delta_\sigma$ be the discrete magnetic Laplacian (DML) 
associated to
a multiplier function $\sigma$ with values in $U(1)$. In particular,
$\sigma$ is a function defined on oriented edges of $K$ and satisfies
$\sigma([u,v])=\overline{\sigma([v,u])}$. Explicitly,
\begin{equation}\label{dml}
\dml f(v) = \sum_{w \sim v} \left (f(v) -
\sigma([w,v])f(w)\right ).
\end{equation}

Our analog of Kato's inequality is as follows.
\begin{lemma} \label{Kato}
For every $f \in C^0_{(2)}(K)$, one has the pointwise inequality:
$$
|f|\cdot \Delta |f| \leq \Re (\dml f \cdot \overline {f}).
$$
\end{lemma}

\begin{proof}
By an explicit calculation using (\ref{lapl}) and (\ref{dml}),
$$
(|f|\cdot \Delta |f|)(v) - \Re (\dml f \cdot \overline {f})(v) =
\sum_{w \sim v} \Re ( \sigma([w,v])f(w)\overline{f(v)} -
|f(v)|\cdot |f(w)|) \leq 0.
$$
\end{proof}

Recall that the bottom of the spectrum of a self-adjoint operator 
such as the DML is computed as:
$$
\lambda_0(\dml) = 
\inf\left\{ \frac{ (\dml f , f)}{||f||^2} : ||f||
 \ne 0\right\}
$$
where $||f||$ denotes the $\ell^2$-norm of $f \in C^0_{(2)}(K)$ and 
$ (\dml f , f)$ the $\ell^2$-inner product. As an immediate 
corollary of Lemma \ref{Kato}, one has,
\begin{corollary}\label{bos}
On $ C^0_{(2)}(K)$, one has the inequality
$$
\lambda_0(\Delta)  \le \lambda_0(\dml).
$$
\end{corollary}

Recall the theorem of Kesten and Brooks \cite{Brooks1982}, which says that for finitely 
generated, infinite groups, one has $\lambda_0(\Delta) = 0$ if and only 
if $G$ is an amenable group. The following is an immediate corollary of this 
fact and Corollary \ref{bos}.

\begin{corollary}
If $G$ is a nonamenable group, then on $ C^0_{(2)}(K)$, one has
$$
 \lambda_0(\dml) > 0.
$$
\end{corollary}

Note that we regard $C^0_{(2)}(K)$ as a Hilbert space with the inner
product $$(f,g) = \sum_v f(v)\cdot \overline {g(v)}.$$ 
Therefore every
bounded operator on it can be represented as a matrix indexed by ordered
pairs of vertices of $K$. Since both Laplacians are bounded operators,
the corresponding heat semigroups $e^{-t\Delta}$ and $e^{-t\dml}$ are
unambiguously defined for $t > 0$. Let $p_t(v,w)$ and $p_t^\sigma(v,w)$
be the corresponding kernels, i.e.\ matrices.

We have the following domination relation.

\begin{theorem} \label{domination}
For every multiplier $\sigma$ and every $f\in C^0_{(2)}(K)$ we have
\begin{equation}\label{heat-domination}
|e^{-t\dml}f| \leq e^{-t\Delta}|f|
\end{equation}
pointwise. As a consequence,
$$
\left | (e^{-t\dml}f,g) \right | \leq (e^{-t\Delta}|f|,|g|)$$
for all $f,g \in C^0_{(2)}(K)$. In particular,
$$
p_t^\sigma(v,v) \leq p_t(v,v)$$
for every vertex $v$ of $K$ and every $t> 0$.
\end{theorem}
\begin{proof}
It follows from the inequality in Lemma \ref{Kato} that for every
$\lambda > 0$
$$
|g|\cdot (\Delta + \lambda I)|g| \leq \Re ((\dml+\lambda I)g\cdot
\overline{g}) \leq | (\dml+\lambda I)g|\cdot |g|$$
for all $g \in C^0_{(2)}(K)$. It follows that
$$
(\Delta + \lambda I)|g| \leq |(\dml + \lambda I)g|.$$
Strictly speaking we can make this conclusion only for vertices where
$g(w)\neq 0$. However, when $g(w)=0$, the left-hand side is nonpositive
by (\ref{lapl}) and the left-hand side is trivially nonnegative.

Now let $g=(\dml + \lambda I)^{-1} f$. The last inequality can be
rewritten as $$
(\Delta + \lambda I)|(\dml + \lambda I)^{-1} f| \leq |f|.$$
Since $(\Delta +\lambda I)$ is positivity preserving we obtain the
(pointwise) inequality
\begin{equation}\label{res-comp0}
|(\dml + \lambda I)^{-1}f| \leq (\Delta + \lambda I)^{-1} |f|.
\end{equation}
By induction,
\begin{equation}\label{res-comp}
|(\dml + \lambda I)^{-n}f| \leq (\Delta + \lambda I)^{-n} |f|
\end{equation}
for every positive interger $n$. Since $e^{-tA} = \lim_{n\rightarrow
\infty} (I+(t/n)A)^{-n}$, we conclude that
\begin{equation}\label{hk}
|e^{-t\dml}f| \leq e^{-t\Delta}|f|.
\end{equation}
This proves the first inequality in the statement of the theorem.
The second inequality follows easily from the definition of the inner
product. Finally, the last inequality is obtained by setting
$g=f=\delta_v$.
\end{proof}

\begin{remark}
It is very easy to see (using the equality $\int_0^\infty
e^{-t(x+\lambda)} dt = (x+\lambda )^{-1}$) that the inequality
(\ref{res-comp0}) used to prove (\ref{heat-domination}) is actually
equivalent to it.
\end{remark}

We next define the Novikov-Shubin type invariants for the DML. We now
assume that the graph $K$ is equipped with a free action of a discrete
group $G$ so that the quotient $K/G$ is finite. If
$\{\,v_1,v_2,\ldots, v_m\,\}$ is a fundamental set of vertices of $K$
and $A$ is a bounded operator on $C^0_{(2)}(K)$ that
is in commutant of the projective $(G, \sigma)$-action on $C^0_{(2)}(K)$
(cf. \cite{MY} for a description of this action), for instance
the DML, then the von Neumann trace of $A$ is given by
\begin{equation}\label{trace}
\tr_{G, \sigma} (A) = \sum_{i=1}^m (A\delta_{v_i},\delta_{v_i} ).
\end{equation}
Let $\theta_{G, \sigma} (t) =\tr_{G, \sigma} (e^{-t\dml}) $ denote the  
von Neumann theta function.
Choosing a weakly equivalent multiplier determines a unitarily 
equivalent operator
and the von Neumann trace of the new operator remains the same, 
therefore
  $\theta_{G, \sigma} (t) =  \theta_{G, [\sigma]} (t)$.
Explicit calculations tend to show that the large time asymptotics of 
$\theta_{G, [\sigma]} (t)$
are of the form $O(t^{-\beta})$ for some $ \beta>0$. This
motivates the following definitions of the {\em Novikov-Shubin type 
invariants}
associated to the DML,
\begin{equation}
\beta(G,[\sigma])=\sup\left\{\beta\in
I\!\!R:\theta_{G, [\sigma]} (t) \;\; \mbox{ is } 
\;\;O(t^{-\beta})\;\mbox{
as }\;t\rightarrow\infty\right\}\in[0,\infty ].
\end{equation}
The special case $\beta(G,0)= \beta(G)$ gives the usual Novikov-Shubin 
invariants
of the Laplacian. Let ${\rm growth}(G)$ denote the growth rate 
(exponent) of balls of large
metric balls in the group with respect to a word metric. More precisely,
$$
{\rm growth}(G) = \lim_{r\to \infty} \frac{\ln(V(r))}{\ln(r)}
$$
where $V(r)$ denotes the volume of a ball of radius $r$ in a word metric.
We assume here that the group $G$ is finitely generated with a
symmetric, finite set of generators $S$ and the distance
of the group element $g$ from the identity is the lenght of the shortest
word in letters from $S$ representing it. The volume is simply the
cardinality.
We recall the following well known theorem of Varopoulos, cf. 
\cite{VCC}.

\begin{theorem}\label{Var}
Let $G$ be a finitely generated discrete group. Then one has
$$\beta(G)
=\frac{{\rm growth}(G)}{2}$$
\end{theorem}

It follows from \eqref{hk} that for every $t >0$ and for every 
$[\sigma]$, one has
\begin{equation}\label{theta}
\theta_{G, [\sigma]} (t) \le \theta_{G, 0} (t).
\end{equation}
Therefore one has:
\begin{corollary}\label{decay}
Let $G$ be a finitely generated discrete group and $\sigma$ a 
multiplier on $G$. Then
$$
\beta(G,[\sigma])
\ge\frac{{\rm growth}(G)}{2} >0.$$
That is, there is a positive constant $C$ independent of $t$ such that 
for all $\; t\gg0$, one has
$$
\theta_{G, [\sigma]} (t) \le C \; t^{-\frac{{\rm growth}(G)}{2} }
$$
\end{corollary}


\subsection{Determinant of the magnetic Laplacian}
It follows therefore from (\ref{res-comp}) that for every $\lambda >0$ 
and every positive
integer $n$
\begin{equation}\label{traces}
\tr_G{(\Delta + \lambda I)^{-n}} \geq \tr_{G, \sigma}{(\dml + \lambda 
I)^{-n}}.
\end{equation}

\begin{proposition} \label{posdet}
Suppose $\ln\det_{G} (\Delta) > -\infty$ then $\ln\det_{G, \sigma}( 
\dml)
  > -\infty$.
\end{proposition}
\begin{proof} Denote by $F(\lambda)$ and $H(\lambda)$ the spectral
density functions of $\Delta$ and $\dml$ respectively. Multiplying both
operators by a suitable constant we can assume that both have norms
bounded by one and their density functions are constant for $\lambda >
1$. We also note that $F(0)=H(0)=0$ since the graph under consideration
is infinite.

Recall that
\begin{equation}\label{detdef}
\ln{\det}_{G} (\Delta) = \int_{0^{+}}^1 \ln \lambda\; dF(\lambda)
\end{equation}
\begin{equation}
\ln{\det}_{G, \sigma} (\dml) =\int_{0^{+}}^1 \ln \lambda\; dH(\lambda).
\end{equation}
We begin by observing that a necessary and sufficient condition for
$\int_{0^{+}}^1 \ln \lambda\; dH(\lambda) >
-\infty$ is that $\int_{0^{+}}^1 (1-\lambda)^{-1}\ln \lambda\;
dH(\lambda)> -\infty$. Now, for $0<x\leq 1$,
\begin{equation}\label{riem-sum}
-\frac{\ln x}{1-x} = \frac{1}{1-x}\int_x^1\frac{du}{u}
      = \lim_{n\rightarrow\infty}\sum_{k=1}^n \frac{1}{n-k}\cdot
      \frac{1}{x+k/(n-k)}.
\end{equation}
The last term in the sum ought to be interpreted as 1/n.
This is simply a Riemann sum approximation of the integral but we remark
that it has been chosen so that the approximation is from below and that
passing to the subsequence $n=2^\ell$ we obtain an approximation of the
function $-(\ln x)/(1-x)$ by a sequence of positive functions that is
increasing in $\ell$. For the remainder of the proof we shall write $n$
for $2^{\ell}$.

Now \begin{multline}\label{lapl-bound}
\infty > -\int_{0^{+}}^1 (1-\lambda)^{-1}\ln \lambda\;
dF(\lambda) = \\
\int_{0^+}^1 \lim_{n\rightarrow\infty}\sum_{k=1}^n \frac{1}{n-k}\cdot
\left (\lambda + \frac{k}{n-k}\right )^{-1}\, dF(\lambda) \\
= \lim_{n\rightarrow\infty} \sum_{k=1}^n \frac{1}{n-k}\cdot \tr_G\,
\left (\Delta + \frac{k}{n-k}I\right )^{-1},
\end{multline}
where we have used the monotone convergence theorem to interchange the
limit and integration and the fact that $F(0)=0$ for the equality
$$
\int_{0^+}^1 \left (\lambda + \frac{k}{n-k}\right )^{-1}\, dF(\lambda)
= \tr_G\,\left (\Delta
+ \frac{k}{n-k}I\right )^{-1}.
$$

By our comparison (\ref{traces}) of the trace of the resolvent for 
$\Delta$ and $\dml$
\begin{multline*}
\sum_{k=1}^n \frac{1}{n-k}\cdot \tr_G\,\left (\Delta +
\frac{k}{n-k}I\right )^{-1} \geq
\sum_{k=1}^n \frac{1}{n-k}\cdot \tr_{G, \sigma}\,\left (\dml + 
\frac{k}{n-k}I\right )^{-1} \\
=\int_{0^+}^1 \sum_{k=1}^n \frac{1}{n-k}\cdot
\left (\lambda + \frac{k}{n-k}\right )^{-1}\, dH(\lambda).
\end{multline*}
The integrands in the last integral form a monotone sequence converging
to \\
$-(1-x)^{-1}\ln x$ and by (\ref{riem-sum}) and (\ref{lapl-bound}) the 
integrals are bounded by $$
-\int_{0^{+}}^1 (1-\lambda)^{-1}\ln \lambda\;
dF(\lambda).$$
Applying the monotone convergence theorem again we see that $$
-\int_{0^{+}}^1 (1-\lambda)^{-1}\ln \lambda\;
dH(\lambda) < \infty,$$
so that $\int_{0^+}^1 \ln \lambda \; dG(\lambda) > -\infty.$

\end{proof}

\begin{proposition}
Let $G$ be an infinite but finitely generated discrete group and 
$\sigma$ a multiplier on $G$. Then
$\ln{\det}_{G, \sigma}( \dml )> -\infty$.
\end{proposition}

\begin{proof}
By the Tauberian theorem in \cite{GS}, we know that Theorem \ref{Var} 
is equivalent to
the following: there are positive constants $C_1, C_2$ such that
\begin{equation}\label{spectraldecay}
C_1\; \lambda^{\frac{{\rm growth}(G)}{2}} \le F(\lambda) \le C_2 \; 
\lambda^{\frac{{\rm growth}(G)}{2}}.
\end{equation}
Integrating equation \eqref{detdef} by parts, one obtains
\begin{equation} \label{three}
\ln{\det}_{G}( \Delta)=   \lim_{\epsilon \rightarrow 0^+}
\Big\{(- \log \epsilon) \big( F (\epsilon)
-F(0) \big) - \int_\epsilon^{1} \frac {F (\lambda)} \lambda d 
\lambda\Big\}
\end{equation}
Using the fact that $  \liminf_{\epsilon \rightarrow 0^+} (- \log 
\epsilon)
\big( F (\epsilon) - F (0) \big) \ge 0$
(in fact, this limit exists and is zero)
and $\frac {F (\lambda)}
\lambda \ge 0$ for $\lambda > 0,$ one sees using equation \eqref{decay} 
that
\begin{equation}\label{four}
\ln{\det}_{G}( \Delta ) \ge -
\int_{0^+}^{1}
\frac {F(\lambda) } \lambda d \lambda \ge -\frac{2C_2}{{{\rm 
growth}(G)}} > -\infty.
\end{equation}
The proof of the Proposition is completed by an application of  
Proposition \ref{posdet}.
\end{proof}


\section{Existence and properties of a ground state}

In this section we follow the notation and formalism of \cite{d2}. An
excellent discussion of analogous questions in the continuous case can
be found in \cite{sullivan-lambda}.
Let $K$ be an infinite, connected graph of bounded valence, i.e.\ $m(x)
\leq M$ for a constant $M>0$ independent of $x\in K$. We abuse the
notation slightly and write $K$ for both the graph and the set of its
vertices. The combinatorial Laplacian $\Delta$ on $K$ or a subgraph of
$K$ acts on functions by the formula (\ref{lapl}).
Our sign convention makes $\Delta$ a nonnegative, self-adjoint operator.
Let $K_1$ be a subgraph of $K$. We denote by $\overset{o}{K_1}$ the set
of vertices of $K_1$ all of whose neighbors are in $K_1$. $\partial K_1$
is the complement of $\overset{o}{K_1}$ in $K_1$.

The following elementary lemmata have been proved in \cite{d2}.

\begin{lemma} \label{max}
Suppose $u(x)$ is a real-valued function on $K_1$ and $x\in \overset{o}{K_1}$.
If $\Delta u(x) < 0$, then $u$ does not have a maximum at $x$.
\end{lemma}

\begin{lemma}\label{harnack}
Suppose $u(x)$ is a function on $K_1$ and $u>0$, $\Delta u > 0$ 
on $\overset{o}{K_1}$. If $x \sim y$ and both $x$ and $y$ are in 
$\overset{o}{K_1}$, then $$
\frac{1}{m(y)} \leq \frac{u(y)}{u(x)} \leq m(x).$$
\end{lemma}
In the sequel, we refer to Lemma \ref{max} and Lemma \ref{harnack}
as the maximum principle and the Harnack inequality respectively.

Using these lemmata we prove the existence of a ground state.

\begin{theorem}\label{ground-state}
Let $\lambda_0=\inf\spec \Delta$. There is a positive eigenfunction 
$\phi$ of $\Delta$ on $K$ belonging to $\lambda_0$, i.e. $$
\Delta \phi = \lambda_0 \phi, \qquad \phi > 0.$$
\end{theorem}
\begin{remark} $\phi$ need not be square-summable. The cone of such
functions need not be one dimensional in
general. We do not know whether $\phi$ is necessarily bounded under our
assumption on $K$.
\end{remark}

\begin{proof}
Consider an exhaustion $K_n$ of $K$ by an increasing sequence of 
finite subgraphs $K_n$ such that $\overset{o}{K_n}\subset K_{n+1}$. 
On each $K_n$
consider the combinatorial Laplacian $\Delta_n$ acting on functions that 
vanish on the boundary of $K_n$. Call $\lambda_n$ the smallest eigenvalue 
of the combinatorial Laplacian $\Delta_n$. By \cite{d2}, $\lambda_n >
\lambda_0 \geq 0$ and  $\lambda_n$ converges to $\lambda_0$. Let $\phi_n$ be the
eigenfunction of $\Delta_n$ belonging to $\lambda_n$ normalized so that 
$0 < \phi_n $ inside $K_n$ and $\phi_n(x_0) = 1$ for some fixed point
$x_0$ that belongs to all $K_n$. Note that $\Delta \phi_n (x) = \lambda_n
\phi_n(x) > 0$ in the interior of $K_n$. 
Since $\phi_n(x_0) =1$, it follows now from the Harnack inequality and 
the connectedness of $K$ that 
the sequence $\phi_n(x)$ is bounded for every fixed $x \in K$. Using
the diagonal process we choose a subsequence, still denoted $\phi_n$ that
converges at \emph{all} points of $K$. Let $\phi$ be the limit of this
subsequence. Since $\lim_{n\rightarrow \infty} \phi_n = \phi$ and 
$ \lim_{n\rightarrow \infty} \lambda_n =\lambda_0$ $$
\Delta \phi = \lambda_0 \phi .$$
The limit function $\phi$ satisfies $\phi \geq 0$, $\phi(x_0)=1$. It is 
strictly
positive by the maximum principle applied to $-\phi$. 
\end{proof}

Note that the ground state $\phi$ constructed above can be identically
equal to one if $\lambda_0 = 0$. However, in the nonamenable case,
$\lambda_0 > 0$ and $\phi$ is nonconstant. 

We denote by $d(x,y)$ the distance between $x$ and $y$ in $K$ i.e.\ the
length of shortest chain of edges connecting $x$ with $y$. The following
lemma is an easy consequence of the Harnack inequality.

\begin{lemma}\label{growth} If $\phi$ is the ground state constructed
above so 
that $\phi(x_0)=1$ then,
for all $x\in K$, $\phi$ satisfies  
$$
M^{-d(x,x_0)} \leq \phi(x) \leq M^{d(x,x_0)},
$$
where $M$ is the uniform bound on the valence of vertices of $K$.
\end{lemma}
From now on we refer to the normalized ground state $\phi$ as above as
the ground state.

To prove a further property of the ground state $\phi$ (so called
completeness) we investigate its behavior under the heat semigroup.
Since the combinatorial Laplacian is a bounded operator on $C^0_2(K)$,
$\|\Delta\| \leq 2 \sqrt{M}$ with our definition of the inner product,
$e^{-t\Delta}$ is defined unambiguously as 
\begin{equation}\label{semigroup}
P_t = e^{-t\Delta} = \sum_{n=0}^\infty\frac{(-t)^n\Delta^n}{n!}.
\end{equation}
Its kernel (matrix) is given by $$
p_t(x,y) = (P_t \delta_x,\delta_y).$$

We need to estimate the behavior of $p_t(x,y)$ for bounded $t$ and large
$d(x,y)$.

\begin{lemma}\label{heat-decay}
For every $T>0$ there exist a constant $C(T)>0$ such that $$
p_t(x,y) \leq \frac{C}{d(x,y)!}$$
for all $t\in [0,T] $.
\end{lemma}
\begin{proof}
Write $\Delta^n(x,y)$ for the matrix coefficient of the $n$-th power of
$\Delta$. Then $\Delta(x,y)=0$ if $d(x,y)>1$ by the definition of the
combinatorial Laplacian. It follows, that
$\Delta^n(x,y)=0$ if $d(x,y) > n$. Now suppose that $d(x,y)=m$. It
follows from (\ref{semigroup}) that 
\begin{equation}\label{series}
p_t(x,y) = \sum_{n=m}^\infty\frac{(-t)^n\Delta^n}{n!}.
\end{equation}
Since the combinatorial Laplacian is bounded, $$
|\Delta^n(x,y)| = |(\Delta^n\delta_x,\delta_y)| \leq c^n.$$
Therefore the series obtain by factoring out $1/m!$ from (\ref{series})
is easily seen to be uniformly bounded for $t\leq T$. This proves the
lemma.
\end{proof}

In view of the decay estimate above, the heat semigroup originally
defined on square-summable functions can be extended to functions of
moderate growth at infinity. Formally, $$
P_tu(x) = \sum_{y\in K} p_t(x,y)u(y).$$
For example, since the number of vertices of $K$ grows exponentially
with distance as does the ground state $\phi$ the series $$
P_t\phi(x) = \sum_{y\in K} p_t(x,y)\phi(y)
$$
converges very rapidly. Our goal is to prove that the ground state
$\phi$ is complete in the sense that \begin{equation}
\label{complete}
e^{-\lambda_0 t} P_t\phi (x) = \phi(x).
\end{equation}
To place this in a proper context, consider the renormalized
semigroup $$\tilde{P}_t=e^{\lambda_0 t}[\phi^{-1}]P_t[\phi],$$ where 
$[u]$ denotes the operator of multiplication by the function $u$.
The infitesimal generator $\frac{d}{dt}|_{t=0}\tilde{P}_t$ is equal to
$-[\phi^{-1}](\Delta - \lambda_0)[\phi]=-L$. A calculation shows that
\begin{equation}\label{generator}
L u (x) = \sum_{y\sim x}\frac{\phi(y)}{\phi(x)} (u(x)-u(y)).
\end{equation}
In order to prove (\ref{complete}) we study bounded solutions of the
initial value problem 
\begin{equation}\label{initial}
\begin{split}
Lu + \frac{\partial u}{\partial t}  &= 0\\
u(x,0) &= u_0(x)
\end{split}
\end{equation}
using the method of \cite{d3}.

We need the parabolic version of the maximum principle.
\begin{lemma}\label{max-par}
Suppose $u(x,t)$ satisfies the inequality $Lu + \frac{\partial
u}{\partial t} < 0$ on $\overset{o}{K}_1 \times [0,T]$. Then the maximum
of $u$ on $K_1\times [0,T]$ is attained on the set $K_1\times \{0\}
\cup \partial K_1\times [0,T]$.
\end{lemma}
\begin{proof}
Suppose $(x_0,t_0)\in \overset{o}{K}_1 \times (0,T]$ is a maximum. It
follows that $\frac{\partial u}{\partial t}(x_0,t_0) \geq 0$ so that 
$Lu(x_0,t_0) < 0$. On the other hand, (\ref{generator}) and  positivity
of $\phi$ imply that $Lu(x_0,t_0) \geq 0$. The contradiction proves the
lemma.
\end{proof}
We now state and prove the main result of this section.
\begin{theorem}\label{exist-unique}
Let $u_0(x)$ be a bounded function on $K$, $|u_0(x)|\leq N_0$. Then the
initial value problem (\ref{initial}) has a unique bounded solution
$u(x,t)$. In addition, $$
|u(x,t)| \leq |N_0|$$
for all $(x,t)$.
\end{theorem}
\begin{proof}
In view of the remarks above about the decay of $p_t(x,y)$ $$
u(x,t) = \tilde{P}_t u_0(x) = \sum_{y\in K} e^{\lambda_0 t} 
p_t(x,y)\frac{\phi(y)}{\phi(x)} u(y)$$
is defined by a rapidly converging series and formal calculations
showing that it satisfies (\ref{initial}) are justified. This gives
existence. Now suppose that $u(x,t)$ is a bounded solution. Let $N_1=
\sup |u(x,t)|$. Fix $x_0 \in K$ and define $r(x) = d(x,x_0)$. By our
assumption on the valence and the Harnack inequality 
\begin{equation} \label{Lr} 
|Lr| < M^2 .
\end{equation}
Consider an auxiliary function $$
v(x,t) = u(x,t) - N_0 - \frac{N_1}{R} (r(x) + M^2t),$$
where $R$ is a large parameter. Let $K_1 = B(x_0,R)$ be the set of
vertices of $K$ at distance at most $R$ from $x_0$. The function
$v(x,t)$ is nonpositive on the set $K_1\times \{0\}
\cup \partial K_1\times [0,T]$ and satisfies 
$(Lu + \frac{\partial}{\partial t}) u < 0$ on $\overset{o}{K}_1\times[0,T]$
because of (\ref{Lr}). Lemma \ref{max-par} implies therefore that
$v(x,t) \leq 0$ so that $$
u(x,t) \leq N_0 + \frac{N_1}{R}(r(x) + Mt)$$
on $B(x_0,R)\times [0,T]$. Keeping $(x,t)$ fixed and letting $R$
increase without bounds, wee see that $
u(x,t)\leq N_0$. Applying the same argument to $-u$ yields $|u(x,t)|
\leq N_0$. Since $T > 0$ and $x$ were arbitrary, this last inequality holds
for all $x\in K$ and $t\geq 0$. Uniqueness follows by considering
the difference of two solutions.
\end{proof}

\begin{corollary} \label{completeness} 
The ground state $\phi$ is complete, i.e.\ $$
e^{\lambda_0 t} \sum_{y\in K}p_t(x,y) \phi(y) = \phi(x)$$
for all $x\in K$ and $t > 0$.
\end{corollary}

\begin{proof}
The function $u(x,t)$ identically equal to one satisfies the initial
value problem (\ref{initial}) with the initial value $u_0(x)$ identically
one. $\tilde{P}_t u_0 (x)$ is another solution with the same initial
value. The two must be equal. Hence $$
e^{\lambda_0 t} \sum_{y\in K} p_t(x,y) \frac{\phi(y)}{\phi(x)} = 1$$
for all $x \in K$ and $t \geq 0$.
\end{proof}

We remark that Lemma \ref{max-par} and Theorem \ref{exist-unique} hold 
for the Laplacian $\Delta$ as well as for the renormalized
Laplacian $L$ since both the maximum principle and the Harnack
inequality hold for $\Delta$. This implies that the fundamental solution
$p_t(x,y)$ for the heat equation on the graph $K$ of bounded valence
is unique since it is a bounded solution of the initial value problem  
\begin{equation}
\begin{split}
\Delta u + \frac{\partial u}{\partial t}  &= 0\\
u(x,0) &= \delta_y.
\end{split}
\end{equation}
This fact is of independent interest and we state it separately as
\begin{theorem} For a graph $K$ of bounded valence, the fundamental
solution of the heat equation is unique. More precisely, the function 
$p_t(x,y)$ of variables $t\geq 0$ and $x,y \in
K$ is the unique bounded, positive solution $v(x,y,t)$ of the equation 
$\Delta_x v
+\frac{\partial v}{\partial t} =0$ satisfying the initial condition 
$v(x,y,0)=\delta_y(x)$. 

\end{theorem}

\section{Long time decay of the heat kernel}

\medskip
We now give a discrete analog of a proof, due to Terry Lyons in 
the continuous case, of the decay of $p_t(x,x)$ as $t\rightarrow \infty$. 
Recall that the graph $K$ is equipped with a free action of a discrete
group $G$ and the quotient $K/G$ is finite. If
$\{\,v_1,v_2,\ldots, v_m\,\}$ is a fundamental set of vertices of $K$
and $A$ is a bounded operator on $\ell^2(K)$ that
is in commutant of the $G$ action on $K$, for instance
the combinatorial Laplacian, then the von Neumann trace of $A$ is given by
$$
\tr_{G} (A) = \sum_{i=1}^m (A\delta_{v_i},\delta_{v_i} ).
$$
Let $\theta_{G} (t) =\tr_{G} (e^{-t\Delta}) $ denote the   
von Neumann theta function.
Define
\begin{equation}
\beta(K)=\sup\left\{\beta\in
I\!\!R:   e^{t\lambda_0(K)} \theta_{G} (t) \;\; \mbox{ is }  
\;\;O(t^{-\beta})\;\mbox{
as }\;t\rightarrow\infty\right\}.
\end{equation}

The continuous analog of the proposition below is contained in the
Appendix to \cite{CCMP}. Our proof is patterned after the argument
there.
\begin{proposition} \label{theta1} 
Under the assumptions above,  $\theta_{G} (t)\leq
Ce^{-\lambda_0 t} t^{-1}$ for large $t$. In particular, 
$\beta(K) \geq 1$.
\end{proposition}

\begin{proof}
If $\lambda_0=\lambda_0(K) = 0$ the statement follows from Theorem
\ref{Var} so we assume that $\lambda_0 > 0$ from now on. Recall that
this implies that the group $G$ is nonamenable.  
By the Theorem \ref{ground-state}, there is a ground state, that is, a 
$\phi>0$ such that $ \Delta \phi =\lambda_{0}
\phi $.  Then either there is a $\gamma\in G$ such 
that
$\gamma^\star \phi$ is not proportional to $\phi$, or for every $\gamma\in G$,
$$ \gamma^\star \phi=\alpha(\gamma)\phi $$ for some morphism
$\alpha: G \rightarrow {\mathbb R}_{+}$.

\noindent{\bf Case 1}.  Suppose that there is a $\gamma\in G$ such that
$\gamma^\star \phi$ is not proportional to $\phi$.  
Then $$ u=\frac{\gamma^\star
\phi}{\phi}>0 $$ is a positive, non-constant harmonic function 
for the Markovian
semi-group with matrix $$ \tilde{p}_{t}(x,y)=e^{\lambda_0
t}p_{t}(x,y)\;\frac{\phi(y)}{\phi(x)}\;. $$ This can be seen as follows.
A simple calculation using the $G$ invariance  $p_t(x,y)=p_t(\gamma x,
\gamma y )$ of the heat kernel yields $$
\tilde{P}_t u(x) = e^{\lambda_0 t} \phi(x)^{-1} \sum_{y\in K}
p_t(\gamma x,z)\phi(z)$$
which, by Lemma \ref{completeness}, is equal to $\frac{\phi(\gamma
x)}{\phi(x)} = u(x)$.

\noindent{\bf Case 2}.  Suppose that for every $\gamma\in G$, $$
\gamma^\star \phi=\alpha(\gamma)\phi $$ for some morphism
$\alpha: G \rightarrow \mathbb R_{+}$.  Then the Markovian semigroup $$
\tilde{p}_{t}(x,y)=e^{\lambda_0 t}p_{t}(x,y)\;
\frac{\phi(y)}{\phi(x)} $$ is
$G$ invariant, by a simple calculation.  Now by \cite{lyons-sullivan}, 
Theorem 3 and the comment at the end of Section 5 of
\cite{lyons-sullivan},
$\tilde{p}_{t}$ admits a non-constant, positive, harmonic function $u$. 

Thus in both cases, $\tilde{p}_{t}(x,y)$ admits a positive, nonconstant,
harmonic function. We use this fact to prove the finiteness of the integral
$\displaystyle \int^{\infty}_{1}e^{\lambda_{0} t}p_{t}(x,x)dt$ as
follows.
Observe that 
$$e^{\lambda_{0}t}p_{t}(x,x)=\left 
(e^{-t(\Delta-\lambda_0)}\delta_x,\delta_x\right ).$$
It follows that the function $t\rightarrow e^{\lambda_{0} t}p_{t}(x,x)$ is
non-increasing. Therefore the convergence of the integral 
$\displaystyle \int^{\infty}_{1}e^{\lambda_{0} t}p_{t}(x,x)dt$ is
equivalent to the convergence of the sum 
$\displaystyle
\sum^{\infty}_{n=1}e^{\lambda_{0} n}p_{n}(x,x)$. Now the restriction of
$\tilde{p}_t(x,y)$ to positive integer values of $t$ defines a Markov chain with
transition probabilities $q_n(x,y)=\tilde{p}_{n}(x,y)$. Clearly, the
functions harmonic for $\tilde{p}_{t}(x,y)$ are harmonic for $q_n(x,y)$.
By \cite[Proposition 6.3, Chapter 6]{KSK} this Markov chain is
transient, i.e.\ $$
\sum_{n=1}^\infty q_n(x,x) =\sum^{\infty}_{n=1}e^{\lambda_{0} n}p_{n}(x,x)
$$
is finite. Thus $\displaystyle \int^{\infty}_{1}e^{\lambda_{0}
t}p_{t}(x,x)dt < C$
and, using the monotonicity again, we obtain $$ 
\frac{t}{2}\,e^{\lambda_{0} t}p_{t}(x,x)\leq
\int^{t}_{t / 2}e^{\lambda_{0} s}p_{s}(x,x)ds \leq C`.$$ Summation over 
$x$ in a fundamental domain yields 
$\theta_G(t) = O(e^{-\lambda_0 t} t^{-1})$.

\end{proof}
Proposition \ref{theta} gives an estimate for the trace of the
heat kernel of the discrete Laplacian. An anologous estimate holds for 
the magnetic Laplacian. More precisely, we have:
\begin{corollary}
For every  multiplier
$\sigma$, $$
e^{\lambda_0(K)t}\tr_{G,\sigma} (e^{-t\dml} ) = O(t^{-1})$$
for large $t$.
In addition, if $T$ is a bounded operator on $C^0_{(2)}(K)$ that 
commutes with the projective $(G, \sigma)$-action, and satisfies 
$T \ge -\lambda_0(K) $, then 
$$
\tr_{G,\sigma} (e^{-t(\dml + T)} ) = O(t^{-1}).
$$
\end{corollary}
\begin{proof}
The first inequality follows from the discrete Kato inequality, 
$$\tr_{G,\sigma} (e^{-t\dml}) \leq 
\tr_G (e^{-t\Delta})$$ and Proposition \ref{theta1}. The second inequality is a
consequence of the first one in view of \cite[Theorem 2.16]{HSU1977}.
\end{proof}



\begin{thebibliography}{9999}

\bibitem{An}
A.~Ancona.
\newblock Th\'eorie du potentiel sur les graphes et les vari\'et\'es.
\newblock In {\em \'Ecole d'\'et\'e de Probabilit\'es de Saint-Flour
  XVIII---1988}, volume 1427 of {\em Lecture Notes in Math.}, pages 1--112.
  Springer, Berlin, 1990.
  
  \bibitem{Bellissard}
J. Bellissard, H. Schulz-Baldes, A. van Elst. 
\newblock The Noncommutative Geometry of the Quantum Hall Effect
\newblock{\em J. Math. Phys.}, {\bf 35}, (1994)  5373-5471. 

\bibitem{Berard1986}
P.~H. B{\'e}rard.
\newblock {\em Spectral geometry: direct and inverse problems}, volume 1207 of
  {\em Lecture Notes in Mathematics}.
\newblock Springer-Verlag, Berlin, 1986.
\newblock With appendixes by G\'erard Besson, and by B\'erard and Marcel
  Berger.
  
\bibitem{Brooks1982}
R. Brooks.
\newblock The fundamental group and the spectrum of the Laplacian.
\newblock{\em Comment. Math. Helv.} {\bf 56}  no. 4, (1981) 581--598.


\bibitem{CCMP}
A.~L. Carey, T. Coulhon, V. Mathai, and J. Phillips.
\newblock Von {N}eumann spectra near the spectral gap.
\newblock {\em Bull. Sci. Math.}, {\bf 122} no. 3 (1998) 203--242.

\bibitem{d2}
J.~Dodziuk.
\newblock Difference equations, isoperimetric inequality and transience of
  certain random walks.
\newblock {\em Trans. Amer. Math. Soc.}, {\bf 284} (1984) 787--794.

\bibitem{d3}
J.~Dodziuk.
\newblock Maximum principle for parabolic inequalities and the heat flow on
  open manifolds.
\newblock {\em Indiana Univ. Math. J.}, {\bf 32}  no.5 (1983) 703--716.

\bibitem{GS}  M. Gromov and M. Shubin. {\em {Von Neumann spectra near zero}} ,
Geom.Func.Anal. {\bf 1} no. 4, (1991) 375-404.

\bibitem{HSU1977}
H.~Hess, R.~Schrader, and D.~A. Uhlenbrock.
\newblock Domination of semigroups and generalization of {K}ato's inequality.
\newblock {\em Duke Math. J.}, {\bf 44} no. 4, (1977) 893--904.

\bibitem{KSK}
J.~G. Kemeny, J.~L. Snell, and A.~W. Knapp.
\newblock {\em Denumerable {M}arkov chains}.
\newblock Springer-Verlag, New York, second edition, 1976.
\newblock With a chapter on Markov random fields, by David Griffeath, Graduate
  Texts in Mathematics, No. 40.

\bibitem{lyons-sullivan}
T.~Lyons and D.~Sullivan.
\newblock Function theory, random paths and covering spaces.
\newblock {\em J. Differential Geometry}, {\bf 18} (1984) 229--323.

\bibitem{MY} 
V.~Mathai and S.~Yates. 
\newblock Approximating spectral invariants of Harper operators on graphs,
\newblock  \textit{J. Functional Analysis,} {\bf 188}, no. 1 (2002) 111-136.

\bibitem{sullivan-lambda}
D.~Sullivan.
\newblock Related aspects of positivity in {R}iemannian geometry.
\newblock {\em J. Differential Geom.}, {\bf 25} no.3 (1987) 327--351.

\bibitem{Sunada}
T. Sunada,  
\newblock A discrete analogue of periodic magnetic Schr\"odinger
operators  
(English summary). 
\newblock Geometry of the spectrum (Seattle, WA, 1993), 283--299, 
Contemp. Math., {\bf 173}, Amer. Math. Soc., Providence, RI, 1994. 

\bibitem{Sunda-Polly} P. W. Sy and T. Sunada, 
\newblock Discrete Schr\"odinger operators on a graph. 
\newblock{\em Nagoya Math. J.} {\bf 125} (1992), 141--150.

\bibitem{Colin} 
Y. C. de Verdi\'ere.
\newblock Spectres de graphes. 
\newblock Cours Sp\'ecialis\'es [Specialized Courses], 4. 
Soci\'et\'eŽ Math\'ematique de France, Paris, 1998.

\bibitem{VCC} 
N. Varopoulos, L. Saloff-Coste and T. Coulhon.
Analysis and
Geometry on Groups. {\em Cambridge University Press, Cambridge, U.K.} (1992).



\end{thebibliography}

\end{document}